\documentclass[11pt]{article}
\usepackage{amssymb}
\usepackage{graphicx}
\usepackage{amsmath}
\usepackage{multirow}
\usepackage{enumerate}
\usepackage[labelfont=bf]{caption}

\DeclareMathOperator{\II}{II}
\DeclareMathOperator{\Hess}{Hess}
\DeclareMathOperator{\Sing}{Sing}
\DeclareMathOperator{\Ind}{Ind}

\newtheorem{theorem}{Theorem}[section]
\newtheorem{proposition}[theorem]{Proposition}
\newtheorem{lemma}[theorem]{Lemma}
\newtheorem{remark}[theorem]{Remark}
\newtheorem{corollary}[theorem]{Corollary}
\newtheorem{definition}[theorem]{Definition}
\newtheorem{examples}[theorem]{Example}

\begin{document}
\title{{\bf On the geometric structure of certain real algebraic surfaces}}
\author{
Miguel Angel Guadarrama-Garc\'\i a\thanks{{\it E-mail address:} 
mikhail.engel@gmail.com \newline
	Instituto de Matem\'aticas, Universidad Nacional Aut\'onoma de M\'exico, 
	\'Area de la Inv. Cient\'\i fica, Circuito Exterior C.U.,	
	Coyoac\'an 04510, M\'exico D.F., M\'exico} \\ Adriana 
	Ortiz-Rodr\'{\i}guez\thanks{{\it E-mail address:} aortiz@matem.unam.mx
	\newline
	Instituto de Matem\'aticas, Universidad Nacional Aut\'onoma de M\'exico, 
	\'Area de la Inv. Cient\'\i fica, Circuito Exterior C.U.,	
	Coyoac\'an 04510, M\'exico D.F., M\'exico}}

\date{}
\maketitle
\begin{abstract}
\noindent In this paper we study the affine geometric structure of the graph 
of a polynomial $f \in \mathbb{R}[x,y]$. We provide certain criteria 
to determine when the parabolic curve is compact and when the unbounded 
component of its complement is hyperbolic or elliptic.
We analyse the extension to the real projective plane of both fields of 
asymptotic lines and the Poincar\'e index at its singular points at infinity.
We exhibit an index formula for the field of asymptotic lines 
involving the number of connected components
of the projective Hessian curve of $f$ and the number of godrons. As an
application of this investigation, we obtain upper bounds, respectively, for the 
number of godrons having an interior tangency and when they have an exterior 
tangency.
\end{abstract}

\noindent {\small{\it Keywords:} parabolic curve, asymptotic fields of lines, real 
algebraic sur\-faces, qua\-dra\-tic differential forms.}

\noindent {\small {\it MS classification:} 53A15, 53A05, 14P05, 14N10, 34K32, 34G20}

\section{Introduction}

There is a well known classification of the points 
 of a smooth surface immersed in the three-dimensional real
affine (projective or Euclidean) space. Any point belongs to one of the 
following types: elliptic, parabolic or hyperbolic. On generic smooth
surfaces, parabolic points appear along a smooth curve (it may be empty) 
called {\it the parabolic curve of the surface}, whose complement is constituted 
by the elliptic and hyperbolic domains. The configuration of these sets, invariant 
under the action of the affine group (or projective group) on 3-space, is the basic 
affine geometric structure 
of the surface. One of the goals in projective and affine differential geometry 
has been the study of this basic geometric structure for smooth and also for
algebraic surfaces, see for example \cite{arn2, ArnldAst, bb, ker, Pnv, segre}. 

In this paper, we focus on the analysis 
of the basic geometric structure of generic algebraic surfaces in $\mathbb{R}^3$ 
that are the graph of a real polynomial $f\in \mathbb{R}[x,y]$.
When the parabolic curve of such a surface $S_f$ is compact, there is one 
unbounded component $C_u$ in the complement of this curve that plays a relevant 
role in the determination of this structure. The class of this component can be 
either elliptic or hyperbolic and when it is known we can specify, in 
the generic case, the class of the 
other connected components that are on the complement of the parabolic curve.

In section 3, we study the distribution of the elliptic and hyperbolic domains: we 
provide in Theorem \ref{fndeterminehcompact} 
conditions on the homogeneous part of the highest degree of $f$ 
that guarantee the parabolic curve is compact and indicate the class of the 
component $C_u$.
At each hyperbolic point, there are two lines tangent to the surface that have a 
contact of order, at least three with the surface. These lines are called {\it 
asymptotic lines}. A parabolic point has exactly one asymptotic line.

When $f$ is a differentiable function defined on the plane $%
\mathbb{R}^{2}
$, it is usual to consider a projection of the geometric structure of $S_f$ 
into the plane. The image of the parabolic curve under such projection is a plane 
curve called {\it Hessian curve of $f$} that is defined by the equation
Hess$f$($x,y$) = 0. The images of the two fields of asymptotic lines  
are described by {\it the second fundamental form of $S_f$, }
$$ \II_f(dx,dy) = f_{xx}(x,y) dx^2 + 2 f_{xy}(x,y) dx dy + f_{yy}(x,y) 
dy^2.$$

\begin{figure}[ht]    
\begin{center}
\includegraphics[width=1.90in]{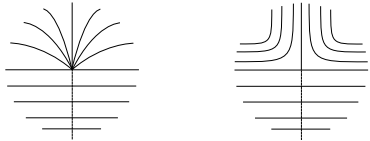}
\caption{Topological types at a singular point at infinity.}
 \label{topologicalpsp}
\end{center}
\end{figure}

\noindent In \cite{Gnnz02}, V. Gu\'\i \~{n}ez considers positive quadratic
differential
equations on the plane $\mathbb{R}^{2}$ of the form
\begin{equation}\label{Pqde}
a\left( x,y\right) dx^{2}+b\left( x,y\right) dxdy+c\left( x,y\right)dy^{2}=0,  
\end{equation}

\noindent where 
$a,b,c\in\mathbb{R}\left[ x,y\right] $ are polynomials of 
degree at most $n$, \,the function  $\,b^{2}-4ac\,$ is nonnegative at every 
point of the $xy$-plane and $b^{2}-4ac \,$ vanishes at a point $p$ if and only 
if $a, b, c$ vanish simultaneously at $p$. 
He extends the fo\-liations determined by equation (\ref{Pqde})
to the line at infinity and proves, among other things, that the 
topological behaviour of these foliations in a neighbourhood of a singular
point at infinity, is one of the types shown in Fig. \ref{topologicalpsp} 
(see \cite{Gnnz02}, 
Remark 2.9).

\noindent When $f\in\mathbb{R}[x,y]
$ is a polynomial, the second fundamental form $\II_f$ 
is a polynomial quadratic differential form that, in general, 
is not positive: 
there are disjoint open sets on the plane where the discriminant of this 
form is negative.

Through the {\it projection of Poincar\'e} from a plane 
into the unitary sphere, we give, in Proposition \ref{extendedqde}, an analytic 
extension on the sphere of the two fields of asymptotic lines. The fields of 
lines obtained through such extension are tangent to the sphere and have the 
same singular points. If the surface $S_f$ is generic, the singular points of
these fields appear on the equator of the sphere and they will be referred to as 
{\it singular points at infinity}. In Theorem \ref{puntos-sing-infinito}, we 
characterise these points. We prove that the Poincar\'e index at a singular point 
at infinity is equal to $\frac{1}{2}$, Theorem \ref{thindexonehalf}. 
As a consequence, we obtain an upper bound for the sum over all Poincar\'e 
indices of an extended field at its singular points. This analysis allowed 
us to itemize the Poincar\'e index at a singular point at infinity
of a field of asymptotic lines when it is extended to the real 
projective plane, Remark \ref{indexatprojplane}.

The projective Hessian curve of $f$ is, in general, a nonsingular algebraic curve 
in $\mathbb{RP}^2$ of even degree. On this plane
we define two surfaces, $B^\pm$, whose boundary is the projective Hessian
curve of $f$. Among parabolic points of a generic surface, $S_f$, a 
{\it godron} 
is distinguished because its unique asymptotic line is tangent to the parabolic 
curve at such point. The problem of determining the lowest upper bound for the 
number of godrons of an algebraic surface in terms of the degree of 
the polynomial that defines it has been an interesting subject of research
\cite{Klkv, salmon, arn1}. The tangency of the asymptotic direction with the 
parabolic curve at such a 
point may be interior or exterior \cite{Bancthom}.
When the surface $B^\pm$ is hyperbolic we give, in Theorem \ref{thglobalbminus}, 
a formula that relates the following three values:
the Euler characteristic of $B^\pm$, the number of godrons having 
either an interior or an exterior tangency and the
Poincar\'e indices at the singular points of the extension to the real 
projective plane of a field of asymptotic lines. 
Derived from this result, upper bounds for the number of 
interior and exterior tangencies are given in Corollary \ref{pssuperiorbounds}.
Another consequence given in Theorem \ref{PSP2014},
is the determination of an upper bound for the number of 
godrons when the projective Hessian curve of $f$ is convex and it is comprised 
only of exterior ovals.
We conclude the paper with the proof of Theorem \ref{thindexonehalf}, section 
\ref{apend}.

\section{Preliminaries}
{\bf\large Classification of points on a generic surface}
\medskip

\noindent A point of a generic smooth surface in $\mathbb{R}^{3}
$ can be classified
in terms of the maximum order of contact of the tangent lines at this point
with the surface \cite{salmon, Lnds, Pltv}. We say that a point $p$ 
is {\it elliptic} if all straight lines tangent to the surface at $p$ have a 
contact of order two with the surface at that point.

An {\it asymptotic line} at a point $p$ is a straight line tangent to the surface 
at $p$ that has a contact of order greater than two with the surface. A {\it 
hyperbolic point} has exactly two transversal asymptotic lines while a {\it 
parabolic point} has one (double) asymptotic line. 

The sets of elliptic and hyperbolic points are open subsets on the surface called 
{\it elliptic} and {\it hyperbolic domains}, respectively. 
These two domains share a common boundary called 
{\it the parabolic curve} which is a smooth curve constituted by the parabolic
points.
The unique asymptotic line at a parabolic point is transversal to the parabolic
curve except at some isolated points called {\it godrons} 
(other authors call them cusps of Gauss or special parabolic points). The order 
of contact of the asymptotic line at 
each parabolic point is three while at a godron, is four. The set of 
asymptotic directions makes up, globally, two continuous fields of directions 
tangent to the surface \cite{Urb} (this property is proved locally in \cite{dara, 
Dvdv}). The integral curves of these fields are known as {\it asymptotic curves}.
A hyperbolic point $p$ is called a {\it point of inflexion} if the order of 
contact of an asymptotic line with the surface at $p$ is at least four. This 
property implies that an asymptotic curve passing through $p$ has an inflexion 
point at such point. The set of points of inflexion is called {\it flecnodal 
curve}. The closure of the flecnodal curve is a curve which is tangent to the 
parabolic curve at the godrons.

In order to understand the geometric structure of $S_f$, when the surface $S_f$ 
is the graph of a differentiable function $f$ on the plane, it is usual to
consider the projection of the elements constituting such structure under the map
$\,\pi : \mathbb{R}^3 \rightarrow \mathbb{R}^2 $, $(x,y,z)\mapsto (x,y)$.
The image of the parabolic curve on the $xy$-plane under
$\pi$ is the zero locus of the Hessian function $\Hess f=f_{xx}f_{yy}-f_{xy}^{2}$.
This curve will be called {\it the Hessian curve of $f$}.
The hyperbolic and elliptic domains are projected, respectively, on 
$H$ and $E$, where
the Hessian function of $f$ is negative and positive, respectively. 
The projection of both fields of asymptotic directions over the $xy$-plane yields 
two fields of lines that are described by the 
quadratic differential equation:
\begin{equation}\label{QDE}
f_{xx}\left( x,y\right) dx^{2}+2f_{xy}\left( x,y\right)
+f_{yy}\left( x,y\right) dy^{2}=0.
\end{equation}
The quadratic form on the left will be referred to as the
{\it second fundamental form of $f$} and will be denoted by $\,\II_f$.
For sake of simplicity, we identify the solutions of this quadratic form
with the asymptotic directions and they will be referred to as {\it the fields of 
asymptotic directions of $f$}. A point on the $xy$-plane is a {\it flat point of 
$\,\II_f$} if the coefficients of this form, $\, f_{xx}, f_{xy}$ and 
$f_{yy}$, vanish at this point. 

We are interested in the particular case when $f\in \mathbb{R}[x,y]$ is a polynomial. 
If the degree of $f$ is $n$, its Hessian curve 
is a real plane algebraic curve of degree, at most $2n-4$. Moreover, if we 
consider {\it the homogeneous decomposition of $f$}, 
$f =\sum_{i=r}^{n}f_{i}$, where $f_i\in \mathbb{R}\left[ x,y\right]$ 
is a homogeneous polynomial of degree $i$, then
$$\mbox{Hess}f(x,y) = \sum_{j=2r-4}^{2n-4} h_{j}\left( x,y\right),
\mbox{ where } h_{2r-4} = \mbox{Hess}f_r \mbox{ and } h_{2n-4} = \mbox{Hess} f_n.$$

\begin{definition}
The projective Hessian curve of $f$ {\rm is the zero locus of the 
homogeneous polynomial $H_f \in \mathbb{R}[x,y,z]$ which is the homogenization 
of the polynomial Hess$f(x,y)$.}
\end{definition}
It follows, from the homogeneous decomposition of $f$, that $H_f$ has the expression:
$\,H_f \left(x,y,z\right) =
\sum_{j=2r-4}^{2n-4} z^{2n-4}\, h_{j}\left(\frac{x}{z},\frac{y}{z}\right)$. 
Therefore, the restriction of $H_f$ to the line at infinity $z=0$ is 
$$H_f(x,y,0) =  \mbox{ Hess} f_n(x,y).$$ 

\noindent When the degree of $H_f$ is even it allows us to label the points at infinity.
A point on the line at infinity $z=0$ of $\mathbb{RP}^{2}$
is called \emph{elliptic, parabolic} or {\it hyperbolic} if the sign of the
homogeneous polynomial $H_f$ at this point is positive, zero or negative,
respectively.

Now, we shall introduce the concept of generic surface for $S_f$ but 
before we give some definitions. 

Let $S$ be a smooth surface in $\mathbb{RP}^3$ and $p$ a point on $S$.
Two  function germs of $S$ at $p$ are {\it equivalent} if one is transformed into the other under the diffeomorphism group action.
In paper \cite{Pltv}, O. A. Platonova proves that \medskip

{\it ``In the space of compact smooth 
surfaces in $\mathbb{RP}^3$ there is an open everywhere dense set of surfaces of 
which the germs at each point are equivalent to the germs that have the $p$-jets 
in Table \ref{platonova}"}.

\begin{table}[ht]
\centering
\begin{tabular}{ | c | c | c | c | c |}\hline
\mbox{Notation} & \mbox{Normal form} & \mbox{Restrictions} & $p$ & \mbox{cod} \\ 
\hline
$\Pi_2$ & $x^2+y^2$ & $-$ & 2& 0 \\ \hline
$\Pi_{3,1}$ & $xy+ s x^3 + y^3$ & $s \neq 0$ & 3 & 0 \\ \hline
$\Pi_{3,2}$ & $y^2 + x^3$ & $-$ & 3 & 1 \\ \hline
$\Pi_{4,1}$ & $xy+ y^3 + x^4 + hx^3y$ & $-$ & 4 & 1 \\ \hline
$\Pi_{4,2}$ & $y^2 + x^2y + \nu x^4$ & $\nu \neq 0, 1/4$ & 4 & 2 \\ \hline
$\Pi_{4,3}$ & $xy + x^4 + s_1 x^3y + s_2 xy^3 + s_3 y^4$ & $s_3 \neq 0$ & 4 & 2 \\ 
\hline
$\Pi_5$ & $ xy + y^3 \pm x^3y \sum d_i x^{5-i} y^i$ & $d_0 \neq 0$ & 5 & 2 \\ \hline
\end{tabular}
\caption{Normal Forms}
\label{platonova}
\end{table}

\begin{definition}{\rm
A smooth surface in $\mathbb{RP}^3$ is {\it generic} if it belongs to the 
open everywhere dense set defined by Platonova.}  
\end{definition}

When $f$ is a polynomial, its graph $S_f$ is an algebraic surface in $\mathbb{R}^3$
and we will say that $S_f$ is {\it generic} if the $p$-jet of the function germ
 at each point of $S_f$ is equivalent to a normal form of Table 
\ref{platonova} and if the projective Hessian curve of $f$ is nonsingular.
\bigskip

\noindent{\bf\large Real algebraic curves in $\mathbb{RP}^{2}$}
\medskip

\noindent A \emph{real algebraic curve} in $%
\mathbb{RP}^{2}$ \emph{of degree $m$} is, up to nonzero constant factors, a 
homogeneous polynomial $%
F\in\mathbb{R}\left[ x,y,z\right] $ of degree $m$. The polynomial equation $F\left(
x,y,z\right) =0$ determines the \emph{set of real points of the curve} in $%
\mathbb{R}
P^{2}.$ From now on, we shall also call this set a real algebraic curve in
$\mathbb{RP}^{2}$.

Each connected component of a nonsingular algebraic curve in $%
\mathbb{RP}^{2}$ is homeomorphic to a circle. There are two ways
up to isotopy to embed a circle into the real projective plane which are called 
the \emph{two-sidedly} and the \emph{one-sidedly} \cite{zeuthen}. 
In the two-sidedly case, the complement in $\mathbb{RP}^{2}$ of the image $L$ 
of the circle has two connected components, 
one of which is homeomorphic to an open disc and called \emph{the inside component 
of $L$}
while the other is homeomorphic to a M\"{o}bius strip and is known as \emph{the 
outside
component of $L$}. Under these conditions, the image of the circle is called 
\emph{oval}.
We say that an oval is an {\it outer oval} if it is not in the inside component of any 
other oval.
In the one-sidedly case, the complement in $\mathbb{RP}^{2}$ of the image of the circle is 
connected and homeomorphic to a disc. 
In this situation, the image of the circle is called \emph{a pseudo-line.} 
While all connected components of a nonempty nonsingular real algebraic 
curve in $\mathbb{RP}^{2}$ of even degree are ovals, each nonsingular algebraic 
curve of odd degree 
is constituted by ovals (if there are any) and exactly one pseudo-line.

\begin{figure}[htb]    
\begin{center}
\includegraphics[scale=.21]{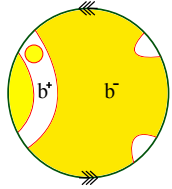}
\caption{Open sets $b^{+}$ and $b^{-}$.}
 \label{complhessproyectiva}
\end{center}
\end{figure}

The complement of a nonsingular curve $F$ of even degree  in $%
\mathbb{RP}^{2}$ is the union of two disjoint open subsets, say $%
b^{+}$ and $b^{-}$ (Fig. \ref{complhessproyectiva}). 
The set $b^{+}$ is an orientable smooth surface at which the sign of $F$ does not 
change while the open set $b^{-}$ is a nonorientable smooth surface at which
$F$ takes the other sign. The closure of $b^{+}$ and $b^{-}$ will be denoted 
by $B^{+}$ and $B^{-}$, respectively.

\begin{definition}{\rm
An oval of a real algebraic curve in $\mathbb{RP}^{2}$ of even degree $m$ 
is called \emph{even} (\emph{odd}) if it is contained 
in an even (odd) number of ovals of the same curve.
The number of even ovals is denoted by $P$ and the number of odd ovals by $N.$ } 
\end{definition}
The numbers $P$ and $N$ contain information about the topology of the
surfaces $B^{+}$ and $B^{-}.$ Indeed, the surface $B^{+}$ has $P$ connected
components and the surface $B^{-}$ has $N+1$ connected components.
In 1906, Virginia Ragsdale proves that the
Euler characteristics of these surfaces are $\chi \left( B^{+}\right) =P-N$
and $\chi \left( B^{-}\right) =N-P+1$, \cite{Rgsdl}. Three decades later, 
I. Petrowsky shows in \cite{Ptrwsk} the following

\begin{theorem}\label{petrows}
Any nonsingular real projective algebraic cur\-ve of even degree $m=2k$
satisfies 
\begin{equation*}
-\frac{3}{2}k\left( k-1\right) \leq P-N\leq \frac{3}{2}k\left( k-1\right) +1.
\end{equation*}
\end{theorem}

\section{Determination of the elliptic and hyperbolic domains}
In this paragraph, we analyse the geometric behaviour of the sets $E$ and $H$.

\begin{definition}{\rm
A homogeneous polynomial on $\mathbb{R}\left[ x,y\right] $\ is called 
$\emph{hyperbolic}$  ($\emph{elliptic}$) 
if its Hessian polynomial has no real linear factors and if it is nonpositive
(nonnegative) at any point.}
\end{definition}

\begin{examples}{\rm \cite{tesismagg}}\label{ejemfactorisables}
{\rm If $f\in \mathbb{R}[x,y]$ is a homogeneous polynomial of degree $n\geq 2$ 
with $n$ real linear factors that are distinct up to nonzero constant factors, 
then it is hyperbolic.}
\end{examples}

Let us denote by $H^n[x, y]\subset\mathbb{R}[x, y]$ the set of real
homogeneous polynomials of degree $n$.
The set constituted by hyperbolic homogeneous polynomials of degree $n$ is a
topological subspace of $H^n[x, y]$, denoted by {\it Hyp}$(n)$. The
connectedness of this space has been studied as part of the subject known
as the Hessian Topology introduced in \cite{arn1, arn2, panov2} and named 
by V. I. Arnold in \cite{ArnldAst}.  
In fact, in reference \cite{ArnldAst} it is shown that
this topological property of {\it Hyp}$(n)$ depends on the degree of the 
polynomials that
constitute it. For example, {\it Hyp}$(3)$ and
{\it Hyp}$(4)$ are connected subspaces whereas {\it Hyp}$(6)$ is a disconnected one.
According to this, V.I. Arnold stated the following conjecture \cite{ArnldAst}, 
p.1067: \medskip

\noindent {\it ``The number of connected components of the space of 
hyperbolic homogeneous polynomials of degree n increases as n increases 
(at least as a linear function of n)."}

\begin{lemma}
$i$) A hyperbolic homogeneous polynomial in $\mathbb{R}[x,y]$ has at least one 
real 
linear factor. Moreover, every real linear factor of a hyperbolic polynomial has 
multiplicity one. \newline
$ii$) An elliptic homogeneous polynomial has no real linear
factors.
\end{lemma}

\noindent {\bf Proof.}
We firstly note that when a real linear factor of a homogeneous polynomial $f$ 
has multiplicity greater than one, 
the Hessian curve of $f$ is unbounded.\newline
$i$) Let suppose that $n$ is even and that $f$ has no real linear factors. 
On the one hand, V.I. Arnold proves in \cite{ArnldAst}
(p.1035) that the index of a field of asymptotic lines at the origin on the $xy$-
plane is
\begin{equation}\label{indice-con-f}
\mbox{ind$_0$(cross($\gamma$))} = 1 - \frac{1}{4} \# \{\theta\in [0, 2\pi) : 
F(\theta) = 0\},
\end{equation}
where $F(\theta)$ is the expression of the hyperbolic homogeneous polynomial in 
polar coordinates and $\gamma$ is a parametrization of the unitary circle 
centred at the origin of the $xy$-plane. 

According to (\ref{indice-con-f}) and considering the fact that 
$f$ has no real linear factors,  
ind$_0$(cross($\gamma$)) = 1.
On the other hand, Arnold shows in the same paper (p.1038) that if $f$ is a 
hyperbolic homogeneous polynomial of degree even, then 
\begin{equation}
\mbox{ ind$_0$(cross($\gamma$))}\leq 0.
\end{equation}
It is a contradiction to the first assertion.

\noindent $ii$)
Let $f$ be an elliptic homogeneous polynomial. In this case, its Hessian curve is 
compact, in fact, it is the origin. Let us suppose that $f$ has a real linear 
factor $l(x,y)$. Thus, $l=0$ is an asymptotic curve
because it has an infinite order of contact with $S_f$ and
the multiplicity of $l$ is one.  This is a
contradiction. \hfill $\Box$

\begin{examples}{\rm (\cite{Olvr2}, p.60)
For each $\mu \in \left( 0,1\right) $ and $\alpha =\pm 1$, the homogeneous 
polynomial $\alpha \left( x^{4}+6\mu x^{2}y^{2}+y^{4}\right) $ is elliptic. 
Therefore, any element of its orbit is elliptic, by considering
the action of $GL\left( 2, \mathbb{R} \right)$ on $H^4[x,y]$.}   
\end{examples}

In Theorem \ref{fndeterminehcompact}, we can appreciate how $f_n$ determines the 
geometric structure of the surface $S_f$ when such homogeneous polynomial is 
hyperbolic or elliptic. In other circumstances, it can be untrue
as shown by the following examples. Consider the polynomials $\,f(x,y)=x^4 +
6x^2y^2 -y^4 +3x^2y -3xy^2 +10y^2 -10x^2\,$ and $\, g(x,y)=x^4 + 6x^2y^2 -y^4 +3x^2y 
-3xy^2 +10y^2 +10x^2$. While they only 
differ by the quadratic homogeneous part, its geometric structure is different 
because in the first case $H$ is contained in $B^-$, and in the second case, $E$ is 
contained in $B^-$. Indeed, the Hessian polynomial of $f$ is
\begin{eqnarray*}
&&\Hess f(x,y)=144x^4 -576x^2y^2-144y^4-72x^3-216x^2y +216xy^2\\
&&-72y^3-36x^2 +36xy+444y^2 +120x+120y-400,
\end{eqnarray*}

\noindent and its restriction to the straight line $x=0$ is a one-variable 
polynomial without real roots. 
On the other hand, the Hessian polynomial of $g$ is
\begin{eqnarray*}
&&\Hess g(x,y)= 144x^4-576x^2y^2-144y^4-72x^3-216x^2y +216xy^2\\
&&-72y^3+228x^2 +36xy+180y^2 -12x+120y+40,
\end{eqnarray*}

\noindent and its restriction to the line $y=0$ is a one-variable 
polynomial without real roots. 
We remark that the intersection of both projective Hessian curves
with the line at infinity are the points $\,P^\pm = [\pm \sqrt{\sqrt{10} -
3}:1:0].$ 
Therefore, the Hessian curve of $f$ has two unbounded connected components: 
one of them is located in both quadrants, 
first and fourth, while the other component is located in the two complementary 
quadrants. This implies that $H$ is contained in $B^-$.
In an analogous way, the Hessian curve of $g$ has two unbounded connected 
components: one of them is located in the first and second quadrants while the 
other is located in the two complementary quadrants. So, the set $E$ is contained 
in $B^-$.

The conclusions of the next result are proved in \cite{H&O&S2} (Theorem 2) 
by considering the
compactification of a plane with the point at infinity. They show that the
fields of asymptotic directions are extended up to the point at infinity
by means of a polynomial binary differential form $\widetilde{\II}_f$. They asume
the extra hypothesis:  the associated form $\widetilde{\II}_f$ (of ${\II}_f$) at
the point at infinity has good multiplicity, that is,  the homogeneous
part of the lowest degree of $\widetilde{\II}_f$ is determined only by the form
${\II}_{f_n}$.

\begin{theorem}\label{fndeterminehcompact}
Let $f \in\mathbb{R}
\left[ x,y\right] $ be a polynomial of degree $n\geq 3$. If $\, f_n$ is 
hyperbolic or elliptic, then the 
Hessian curve of $f$ is compact. Moreover, the set $b^{-}\cap\mathbb{R}^{2}$ 
is hyperbolic or elliptic providing that $\, f_n\,$ is hyperbolic or 
elliptic, respectively. 
\end{theorem} 

\noindent {\bf Proof.}
Suppose that $f_n$ is a hyperbolic polynomial. The elliptic case is similar.
Since the polynomial $\,\Hess f_{n}\left( x,y\right)$ has no real linear factors 
and $\,h_{2n-4} = \Hess f_{n},\,$
the projective Hessian curve does not intersect to the line at infinity.
Accordingly, the curve $\Hess f\left( x,y\right) =0$ is compact in $%
\mathbb{R}^{2}$ and the line at infinity is contained in $B^{-}$.
To show that the set $b^{-}\cap
\mathbb{R}^{2}$ is hyperbolic, it will be enough to prove that
any point on the line at infinity is a hyperbolic point. By taking 
$\, p=\left[ 1:0:0\right]$ we have that $\,H_f(p)$ is negative because 
$\,H_f(p)= \Hess f_{n}\left( 1,0\right)$.
\hfill $\Box$

\section{Projection into the Poincar\'e sphere}
A good approach to studying the behaviour of the asymptotic curves of 
(\ref{QDE}) ``at infinity" is to use the so-called {\it Poincar\'e sphere} \cite{poinc}. 
Let $\,\mathbb{S}^{2} = \{(u,v,w)\in \mathbb{R}^3 \, |\, u^2+v^2+w^2=1\}$ be the unit sphere
centred at the origin $O$ in $\mathbb R^3$ and identify
its tangent plane $T_N\mathbb S^2$ at the north pole $N=(0,0,1)$ with
the $xy$-plane. Given a point ${\bf x}=(x,y,1)\in T_N\mathbb S^2$, the
line through ${\bf x}$ and $O$ intersects $\mathbb S^2$ at the
following two points:

\begin{equation*}
s_{1}\left( \mathbf{x}\right) =\frac{\mathbf{x}}{\sqrt{1+x^{2}+y^{2}}
}, \,\,\,\, s_{2}\left( \mathbf{x}\right)
=-\frac{\mathbf{x}}{\sqrt{1+x^{2}+y^{2}}}.
\end{equation*}%
\begin{equation*}
\end{equation*}
The maps $s_i:\mathbb R^2\to \mathbb S^2$, $i=1,2,$ are called {\it the
projections to the Poincar\'e sphere}.

Now, suppose that $f\in\mathbb{R}\left[ x,y\right] $ is a polynomial of degree $n$
and consider on the $xy$-plane the two 
fields of asymptotic directions, $\mathbb{X}_{1}$ and $\mathbb{X}_{2}$, defined by 
equation (\ref{QDE}).

\begin{remark}
The images of the two fields $\mathbb{X}_{1}$ and $\mathbb{X}_{2}$, under the Poincar\'e projection, over both upper and lower hemispheres,  
are the zero loci of the induced quadratic differential forms, 
$s_{1}^{\ast }\left( \II_f \right)$
and $s_{2}^{\ast }\left( \II_f \right)$, which are defined on the complement of 
the equator of $\,\mathbb{S}^{2}$. Moreover, the images of both fields over 
each open hemisphere
consist of two fields of lines diffeomorphic to $\mathbb{X}_{1}$ and $\mathbb{X}_{2}$. 
\end{remark}

Similarly, as V. Gu\'\i \~{n}ez does in \cite{Gnnz02}, we shall 
prove that the induced quadratic differential forms $s_{1}^{\ast }\left( \II_f
\right)$ and $s_{2}^{\ast }\left( \II_f \right)$ can be extended to an analytical
quadratic differential form defined on the sphere.

\begin{proposition}
\label{extendedqde}
The induced differential forms
$s_{1}^{\ast }\left( \II_f \right) $ and $s_{2}^{\ast }\left( \II_f \right) $ 
are extended
to this analytical quadratic differential form  
\begin{equation}\label{EDLA}
\begin{pmatrix}
du & dv & d\omega%
\end{pmatrix}%
\begin{pmatrix}
\omega ^{2} F_{uu}\left( u,v,\omega \right) & \omega
^{2}F_{uv}\left( u,v,\omega \right) & \omega A\left(
u,v,\omega \right) \\
\omega ^{2}F_{uv}\left( u,v,\omega \right) & \omega
^{2} F_{vv}\left( u,v,\omega \right) & \omega B\left(
u,v,\omega \right) \\
\omega A\left( u,v,\omega \right) & \omega B\left( u,v,\omega \right) & 
S\left( u,v,\omega
\right)%
\end{pmatrix}%
\begin{pmatrix}
du \\
dv \\
d\omega%
\end{pmatrix}%
\end{equation}
\noindent defined on the sphere with the property that the equator is an 
integral curve of the fields defined by this form. 
In such case
\begin{eqnarray*} 
&&F(u,v,\omega) = \sum_{i=0}^{n}\omega ^{n-i} f_i(u,v), \,\, 
F_{uu}= \frac{\partial^2 F}{\partial u^2}, \,\, F_{uv}= \frac{\partial^2 F}
{\partial u\partial v}, \,\, F_{vv}= \frac{\partial^2 F}{\partial v^2},\qquad\qquad\\
&&A\left( u,v,\omega \right) =-u F_{uu}(u,v,\omega)-v F_{uv}
(u,v,\omega), \\
&&B( u,v,\omega) = -u F_{uv}(u,v,\omega) -v F_{vv}(u,v,\omega),\\
&&S\left( u,v,\omega \right) = u^{2} F_{uu}(u,v,\omega) +2uv F_{uv}(u,v,\omega) 
+ v^{2} F_{vv}(u,v,\omega).
\end{eqnarray*}
\end{proposition}

We denote by $\mathbb{Y}_{1},\,\mathbb{Y}_{2}$ the two fields of lines defined
by the form (\ref{EDLA}).
It is worth mentioning that the fields $\mathbb{Y}_{k}, k=1,2,$ are not defined, 
in general, on the whole sphere.
\bigskip

\noindent {\bf Proof.} 
Consider the map $\varrho : \mathbb{R}
^{3}\backslash \left\{ \omega =0\right\} \rightarrow\mathbb{R}
^{2},\,$ $\left( u,v,\omega \right) \mapsto \left(x,y\right)\,$ where
$\, x = \frac{u}{\omega },\, y = \frac{v}{\omega }.$  The images under this map 
of a pair of antipodal points on the sphere $\mathbb{S}^{2}$ are the same.
We proceed to obtain the pullback $\varrho
^{\ast }\left( \II_f \right) $ of the second fundamental form $\II_f .$
Replacing

$\begin{pmatrix}
dx & dy%
\end{pmatrix}%
= 
\begin{pmatrix}
\frac{\omega du-ud\omega }{\omega ^{2}} &
\frac{\omega dv-vd\omega }{\omega ^{2}}%
\end{pmatrix}%
= 
\frac{1}{\omega ^{2}}%
\begin{pmatrix}
du & dv & d\omega%
\end{pmatrix}%
\begin{pmatrix}
\omega & 0 \\
0 & \omega \\
-u & -v%
\end{pmatrix}%
$

\noindent in the expression 
$\, \II_f (dx, dy) = \begin{pmatrix}
dx & dy%
\end{pmatrix}%
\begin{pmatrix}
f_{xx}\left( x,y\right) & f_{xx}\left( x,y\right) \\
f_{xy}\left( x,y\right) & f_{yy}\left( x,y\right)%
\end{pmatrix}%
\begin{pmatrix}
dx \\
dy%
\end{pmatrix}%
,\,$ we have that $\varrho ^{\ast }\left( \II_f \right) $ is
\begin{equation*}
\frac{1}{\omega ^{4}}%
\begin{pmatrix}
du & dv & d\omega%
\end{pmatrix}%
\begin{pmatrix}
\omega & 0 \\
0 & \omega \\
-u & -v%
\end{pmatrix}%
\begin{pmatrix}
f_{xx}\left( \frac{u}{\omega },\frac{v}{\omega }\right) & f_{xy}\left( \frac{%
u}{\omega },\frac{v}{\omega }\right) \\
f_{xy}\left( \frac{u}{\omega },\frac{v}{\omega }\right) & f_{yy}\left( \frac{%
u}{\omega },\frac{v}{\omega }\right)%
\end{pmatrix}%
\begin{pmatrix}
\omega & 0 & -u \\
0 & \omega & -v%
\end{pmatrix}%
\begin{pmatrix}
du \\
dv \\
d\omega%
\end{pmatrix}.
\end{equation*}%
After multiplication by $\,\omega ^{n+2}\,$ we obtain the desired quadratic form
(\ref{EDLA}).

Now, we shall prove that the equator is an integral curve of the fields 
$\mathbb{Y}_{k}, k=1,2$. 
Consider the chart $u=1$ in $\mathbb{R}^3$. In this chart, 
the projections of $\mathbb{Y}_{1}$ and
$\mathbb{Y}_{2}$ restricted to the set $\{(u,v,\omega)\in \mathbb{S}^2 | u> 0\}
$ are described by the quadratic equation

\begin{equation}\label{cartau=1}
\begin{pmatrix}
dv & d\omega%
\end{pmatrix}%
\begin{pmatrix}
\omega ^{2}F_{vv}\left( 1,v,\omega \right) & \omega B\left(
1,v,\omega \right) \\ 
\omega B\left( 1,v,\omega \right) & S\left( 1,v,\omega \right)%
\end{pmatrix}%
\begin{pmatrix}
dv \\ 
d\omega%
\end{pmatrix}%
=0.  
\end{equation}

If the origin of the $v\omega$-plane is not a solution of $S(1,v,\omega)$, the 
following two vector fields are tangent to $\mathbb{Y}_{1}$ and
$\mathbb{Y}_{2}$ in a neighbourhood of the origin,
$$ \frac{d\omega}{dv} = \frac{- \omega B \pm \sqrt{\omega^2 (B^2- F_{vv} S)}}
{S}.$$
So, the $v$-axis is locally an integral curve of $\mathbb{Y}_{k}, k=1,2$. 
\hfill $\Box$ 

\begin{lemma}\label{expresion-S}
The polynomial $S$ of (\ref{EDLA}) is equal to the expression $\, S(u,v,\omega) 
=\sum_{k=2}^{n} k(k-1)\, \omega^{n-k} f_k(u,v).$ 
\end{lemma}

\noindent {\bf Proof.} 
By definition $\, F(u,v,\omega) = \sum_{k=0}^{n}\omega ^{n-k} f_k (u,v)$. Thus,
$$S(u,v,\omega) = 
\sum_{k=0}^{n}\omega ^{n-k} \left(u^2\frac{\partial^2}{\partial u^2}f_k(u,v)
+ 2 uv \frac{\partial^2}{\partial u \partial v}f_k(u,v) + v^2 \frac{\partial^2}
{\partial v^2}f_k(u,v) \right).$$

\noindent By considering the well known Euler's 
formula for a homogeneous polynomial $P \in \mathbb{R}
\left[ x,y\right]$ of degree $m$: $\, m P(x,y)= xP_x(x,y) + yP_y(x,y)$, 
it follows the relation
\begin{equation*}
	m\left( m-1\right) P\left( x,y\right)
	=x^{2}P_{xx}\left( x,y\right)+2xyP_{xy}\left( x,y\right)+y^{2}P_{yy}\left( 
	x,y\right).
\end{equation*}

\noindent We obtain the desired equality by taking $P = f_k$ and $\, m= k$.
\hfill $\Box$

\begin{definition}{\rm
A singular point of $\mathbb{Y}_{k}$ is called {\it singular point at
infinity} if it is on the equator of $\mathbb{S}^{2}$.}
\end{definition}

We remark that if $S_f$ is generic, every singular point of
$\mathbb{Y}_{k}$ is a singular point at infinity.
We say that a point $p \in \mathbb{S}^{2}$ is a {\it flat point of the
quadratic form {\rm (\ref{EDLA})}} if the coefficients of this form vanish
at this point. 

\begin{remark}
(i) A point $\left( u_{0},v_{0},\omega _{0}\right)\in \mathbb{S}^{2}$,
with $\,\omega _{0}\neq 0$ is a flat point of {\rm (\ref{EDLA})} if and
only if the point $\left( x_{0},y_{0}\right) =\left( \frac{u_{0}}{\omega
_{0}},\frac{v_{0}}{\omega _{0}}\right) $ is a flat point of the 
fundamental form $\II_f$.\newline
(ii) By Lemma \ref{expresion-S}: a flat point of {\rm (\ref{EDLA})} 
lies in the equator if and only if the polynomial
$\, S(u,v,0) = k (k-1) f_n(u,v) \,$ vanishes at that point. Therefore, the form
{\rm (\ref{EDLA})} has a finite number of flat points on the equator.
\end{remark}

\begin{theorem}\label{puntos-sing-infinito}
Let $f \in\mathbb{R}\left[x,y\right] $ be a polynomial of degree $n\geq 3$. 
If $p$ is a point on the equator of $ \mathbb{S}^2$, then $p$ is a flat point of 
{\rm (\ref{EDLA})} 
if and only if $p$ is a singular point at infinity of  $\mathbb{Y}_k, \,k=1,2$.
Moreover, if $p$ is a singular point at infinity of  $\mathbb{Y}_k,\,$ and  $f_n$ 
has no repeated factors, then $H_f(p) < 0$.
\end{theorem}

\noindent{\bf Proof}.
Let us suppose that $p = (1,0,0)$. By taking the chart 
$u=1$, the fields $\mathbb{Y}_{k},\, k=1,2,$ restricted to the set $\{(u,v,
\omega)\in \mathbb{S}^2 | u> 0\}$, are described by the quadratic equation
\begin{equation}\label{forma-carta-afin}
\omega^2 F_{vv} (1,v,\omega) dv^2 + 2 \omega B(1,v,\omega) dv d\omega  + 
S(1,v,\omega) d\omega^2 = 0.
\end{equation}

\noindent The discriminant of the left-side form of (\ref{forma-carta-afin}) is 
$\Delta = - \omega^2 \left(F_{vv} S - 
B^2\right)|_{( 1,v,\omega )}$. A straightforward calculation shows that
$\Delta = - \omega^2 H_f(1,v,\omega).$  By Proposition 
\ref{extendedqde}, the $v$-axis is an integral curve of $\mathbb{Y}_k$.
Moreover, since $S(1,v,0)$ has a finite number of solutions, the fields 
of (\ref{forma-carta-afin}) are described, in a neighbourhood of the origin, 
by 
\begin{eqnarray}\label{campos-separados}
R_{k}\left( v,\omega \right) dv + 2 \, S\left(1, v,\omega \right) d\omega =
0,\,\quad k=1,2,
\end{eqnarray}
where $\,
R_{k}\left( v,\omega \right) =- 2 \,\omega B\left( 1,v,\omega \right) + 
2 \left( -1\right)^{k} 
\sqrt{-\omega ^{2} H_f(1,v,\omega )}.$
Suppose that $p$ is a flat point of the form (\ref{EDLA}). So, 
the origin of the $v\omega$-plane is a singular point of the fields  
(\ref{campos-separados}) since it is a zero of $\, 
S(1,v,\omega)$. 
Conversely, let us suppose that the origin of the $v\omega$-plane 
is a singular point of the fields defined by (\ref{campos-separados}). Thus, 
$S(1,0,0) = 0$, that is, the point $p$ is a flat point of 
the form (\ref{EDLA}).  

In order to prove the second part, we suppose that $S(1,0,0) =0$. 
By Lemma \ref{expresion-S}, 
$S(1,0,0) = f_n(1,0) = 0$. It implies that the 
polynomial $y$ is a factor of $f_n$. Thus, by hypothesis the multiplicity of
$y$ is one. In this case $\, H_f(1,y,0) = - (g_x(1,y))^2$, 
where $f_n(x,y) = y g(x,y)$. Thus, $H_f(p) < 0$. Moreover, the discriminant 
$\Delta$ is locally 
positive in the complement of the $\omega$-axis. 
\hfill $\Box$
\medskip

In the next result, we prove that the number of singular points at infinity
of the field  $\mathbb{Y}_{k}, \, k=1,2,$ is twice the number of distinct real
linear factors of the homogeneous polynomial $f_{n}$. Its proof follows from 
Lemma \ref{expresion-S} and Theorem \ref{puntos-sing-infinito}.

\begin{corollary}\label{lemma01}
Let $f \in\mathbb{R}\left[ x,y\right] $ be a polynomial of degree $n\geq 3$. 
Then, the set of singular points at infinity  of $\mathbb{Y}_k,\, k=1,2,\,$ is 
	\begin{equation*}
	\{(u,v,0) \in \mathbb{S}^{2} \, |\, f_{n}\left( u,v\right) =0\}.
	\end{equation*}
\end{corollary} 

\noindent The singular points at infinity of $\mathbb{Y}_k$ that do not belong to the boundary of
$B^\pm$ are characterised in the following

\begin{figure}[htb]    
\begin{center}
\includegraphics[width=0.75in]{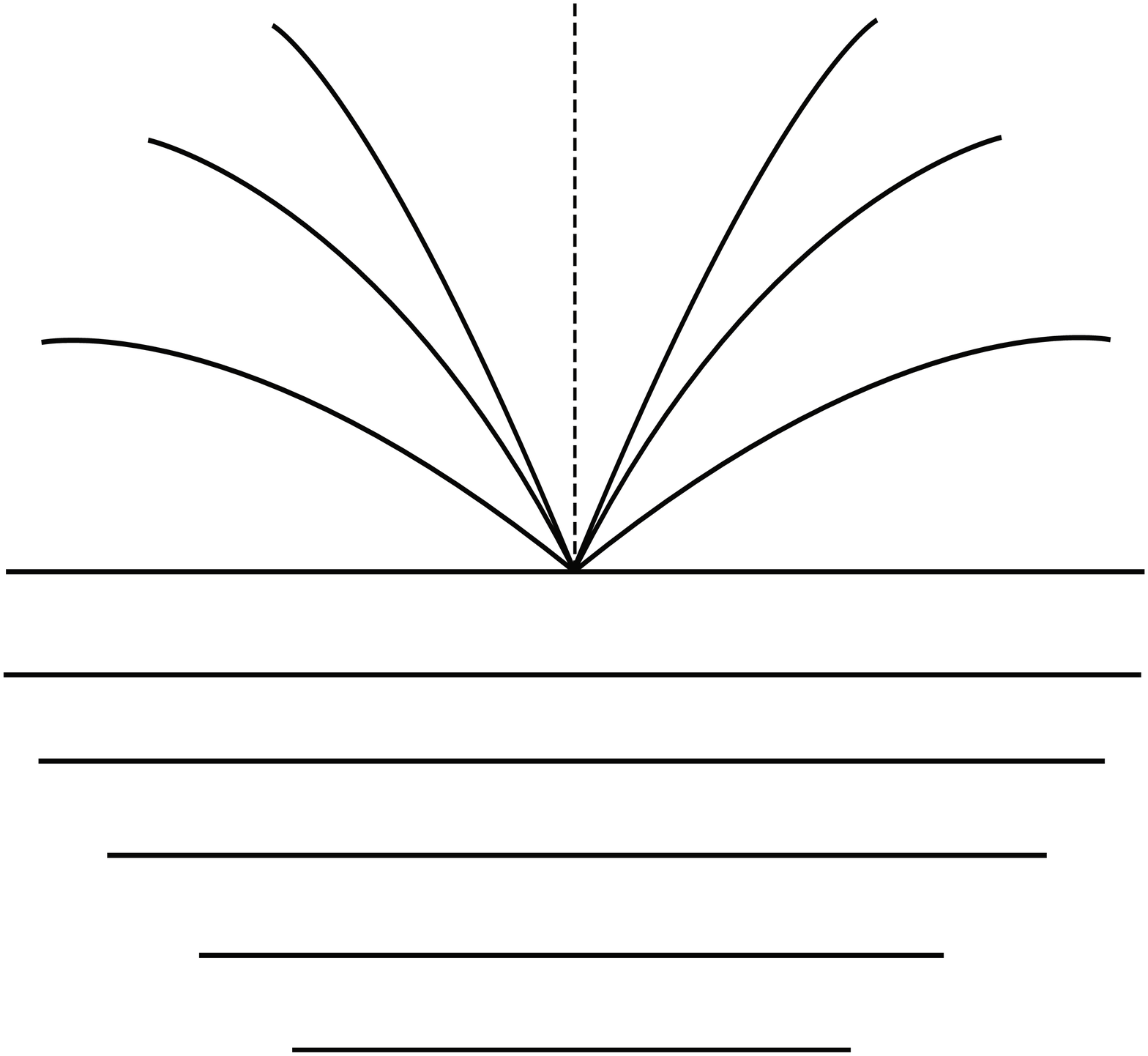}
\caption{Topological behaviour of $\mathbb{Y}_k$ at a singular point at infinity.}
 \label{thetopotype}
\end{center}
\end{figure}

\begin{theorem}\label{thindexonehalf}
Let $f \in \mathbb{R} \left[ x,y\right] $ be a polynomial of degree $n\geq 3$  
such that its homogeneous part $\, f_n$ has no repeated factors.
Then, the Poincar\'{e} index of $\,\mathbb{Y}_k, k=1,2,\,$
at a singular point at infinity is equal to $\frac{1}{2}$. 
Moreover, its topological 
type is the one shown in Fig. \ref{thetopotype}. 
\end{theorem}

\noindent The proof of Theorem \ref{thindexonehalf} is given in section \ref{apend}.

For $k=1,2,\,$ the expression Sing($\mathbb{Y}_{k}$) denotes the 
set of singular points of the field
$\mathbb{Y}_{k}$.

\begin{corollary}\label{cotacasonoacotado}
Let $f \in \mathbb{R}\left[ x,y\right] \,$ be a polynomial such that $f_n$ has 
no repeated factors. \newline
$i$) If the projective Hessian curve of $f$ has a nonempty 
transversal intersection with the line at infinity, then
\begin{equation*}
0 \leq \sum_{\xi \in \Sing (\mathbb{Y}_{k})} \Ind_\xi \left(\mathbb{Y}_{k} 
\right) \leq
n-2, \,\,\,\, \mbox{for} \,\, k=1,2.
\end{equation*} \newline
$ii$) If the projective Hessian curve of $f$ does not intersect to
the line at infinity, then
\begin{equation*}
0 \leq \sum_{\xi \in \Sing (\mathbb{Y}_{k})} \Ind_\xi\left(\mathbb{Y}_{k}\right) 
\leq n, \,\,\,\, \mbox{for} \,\, k=1,2. 
\end{equation*}
\end{corollary}

\noindent{\bf Proof.}
On the one hand, when the polynomial $f_n$ has exactly $n$ generic real 
linear factors, it is hyperbolic (Example \ref{ejemfactorisables}). 
In such a case, the Hessian curve of $f$ is compact by Theorem
\ref{fndeterminehcompact}, and the field $\mathbb{Y}_{k}$ has $2n$ 
singular points at infinity. In this case, in according to Theorem \ref{thindexonehalf}, 
the field $\mathbb{Y}_{k}$ reaches the upper bound of $ii$).
On the other hand, if the Hessian curve of $f$ is unbounded, $f_{n}$ has at most 
$n-2$ real linear factors, and by Lemma \ref{lemma01}, the maximum number 
of singular points at infinity is $2 (n-2)$. 
Inequalities of $i$) and $ii$) follow from Theorem \ref{thindexonehalf}.
\hfill $\Box$
\medskip

\begin{figure}[htb]    
\begin{center}
\includegraphics[width=2.0in]{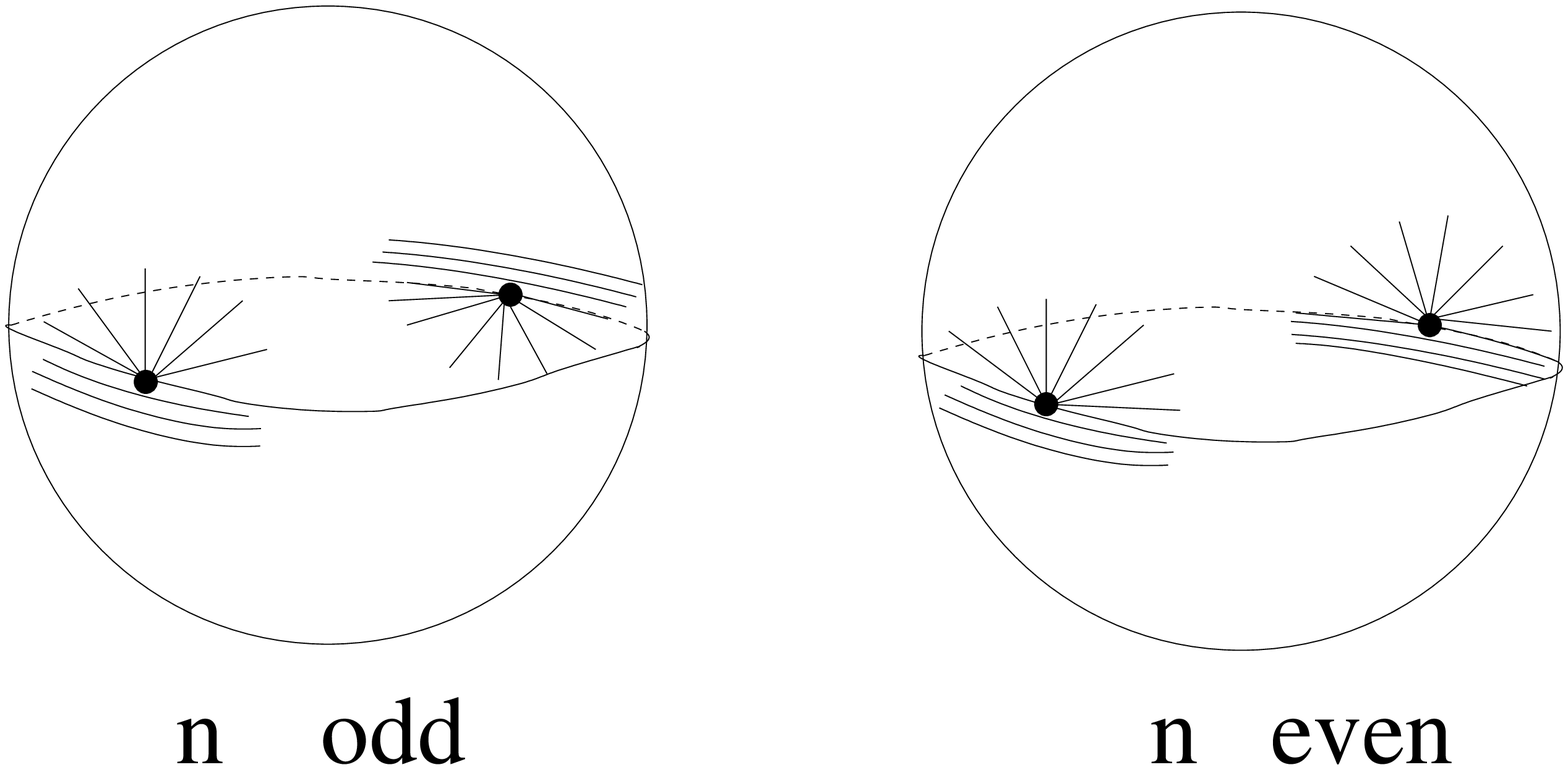}
\caption{Behaviour of $\mathbb{Y}_{k}$ at antipodal singular 
  points at infinity.}
  \label{psingesfera} 
\end{center}
\end{figure}

\begin{remark}\label{foliacionesenesfera}
The fields $\mathbb{Y}_{1}$ and $\mathbb{Y}_{2}$ behave qualitatively as: 
\begin{itemize}
\item When $n$ is even, if $\,\mathbb{Y}_{1}$ is the projection on the upper
hemisphere of $\, \mathbb{X}_{1}$, then, on the lower hemisphere, $\mathbb{Y}_{1}$ is
the projection of $\,\mathbb{X}_{2}$. Analogously for $\mathbb{Y}_{2}$. See Fig.
\ref{psingesfera}.
\item If $\,n$ is odd, the field $\,\mathbb{Y}_{i}$ is the projection over both
hemispheres of either $\mathbb{X}_{1}$ or $\mathbb{X}_{2}$ {\rm (}Fig.
\ref{psingesfera}{\rm )}.
\end{itemize}
\end{remark}

\begin{figure}[htb]    
\begin{center}
\includegraphics[width=1.8in]{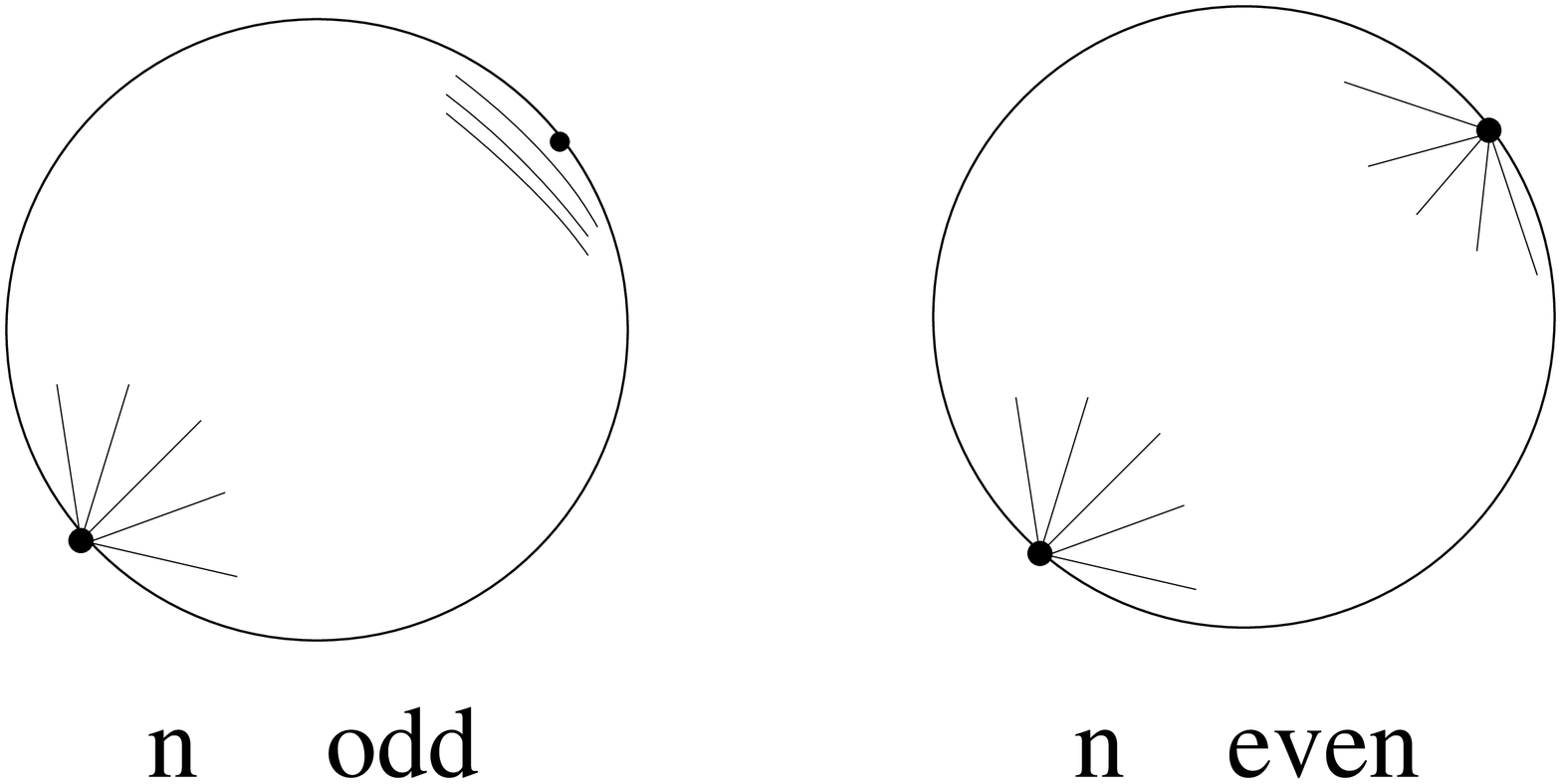}
\caption{Behaviour of $\widetilde{\mathbb{X}_k}$ at a singular point at
     infinity.}
  \label{psingproyec} 
\end{center}
\end{figure}

\noindent The restriction of the field $\mathbb{Y}_{k}, \, k\in\{1,2\},$
to the closure of a hemisphere of $\mathbb{S}^{2}$ will be called {\it a projective 
extension of the field of asymptotic directions} $\mathbb{X}_k$ and it will be denoted 
by $\widetilde{\mathbb{X}_k}$. Let us suppose that $p$ is a singular point at infinity 
of $\mathbb{Y}_k$. When $n$ is odd, a picture of the local qualitative behaviour of any 
projective extension at the points, 
$p$ and $-p$, is shown in Fig. \ref{psingproyec}. We will say that [$p$]$=\{p,-p\}$
is a {\it singular point at infinity of the projective extension}.
Now, let $n$ be even and choose a projective extension $\widetilde{\mathbb{X}_k}$. 
If the local qualitative behaviour of $\widetilde{\mathbb{X}_k}$ at points, $p$ and $-p$, 
is as shown in Fig. \ref{psingproyec}, we will say that [$p$]$=\{p,-p\}$
is a {\it singular point at infinity of $\widetilde{\mathbb{X}_k}$}.
Thus, from Remark \ref{foliacionesenesfera} we have the following

\begin{remark}\label{indexatprojplane}
The Poincar\'e index of a projective extension $\widetilde{\mathbb{X}_k}$ at a 
singular point at infinity is equal to $\frac{1}{2}$ if  
$n$ is odd, and it is 1 when $n$ is even. Their topological types are shown 
in Fig. \ref{psingproyec}.
\end{remark}	

\begin{figure}[htb]    
\begin{center}
\includegraphics[width=2.0in]{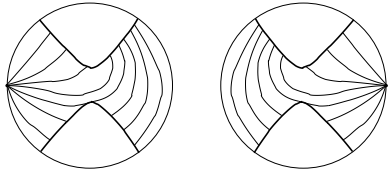}
\caption{Foliation by asymptotic curves for the cubic polynomial $q(x,y)$.}
  \label{ejemplog3} 
\end{center}
\end{figure}

\begin{examples}{\rm 
Consider the cubic polynomial $q(x,y) = x^2 + y^2 + y(x^2 + y^2)$.
In \cite{H&O&S1}, it is proved that the Hessian curve of $q$ is a hyperbola,
it contains one godron and the flecnodal curve of $q$ is the straight line
$y=0$. Moreover, the convex domain is elliptic while the concave is hyperbolic.
By Corollary \ref{lemma01}, each field $\mathbb{Y}_{k}$ has two singular points 
at infinity. In Fig. \ref{ejemplog3}, we draw the foliation of 
$\mathbb{Y}_{k}$ in both closed hemispheres. We remark that this qualitative 
behaviour is the same for any nonhomogeneous cubic polynomial such that $q_3$ 
has exactly one real linear factor.}
\end{examples}

\begin{examples}\label{exfact}{\rm 
The product $\Pi_{k=1}^n l_k$ of $n$ linear polynomials on
$\mathbb{R}[x,y]$ is called {\it a factorisable polynomial of degree $n$} if
($i$) the intersection of each pair of straight lines $l_i=0, \, l_j=0, \,
i\neq j, i,j\in \{1,\ldots ,n\}$ is nonempty and
($ii$) for each $i=1,\ldots ,n,$ the straight line $l_i=0$ has no critical
points of the function $\,\Pi_{j\neq i} l_j$.

The geometrical structure of a factorisable polynomial of degree $n\geq 3$
is described in Theorem 1 of \cite{adriana}. When $n=4$, such geometrical
composition is as follows: the parabolic curve of $f$ is a
quartic smooth and compact curve with three connected components. The
unbounded component $C_u$ is hyperbolic and the graph of $f$ has eight
godrons, all of index -1. Moreover, the flecnodal curve is only
constituted by the straight lines $l_i(x,y)=0, i=1,\ldots ,4$.

\begin{figure}[ht]  
\centering
  \begin{minipage}{0.46\textwidth}
  \quad\qquad\includegraphics[width=0.7in]{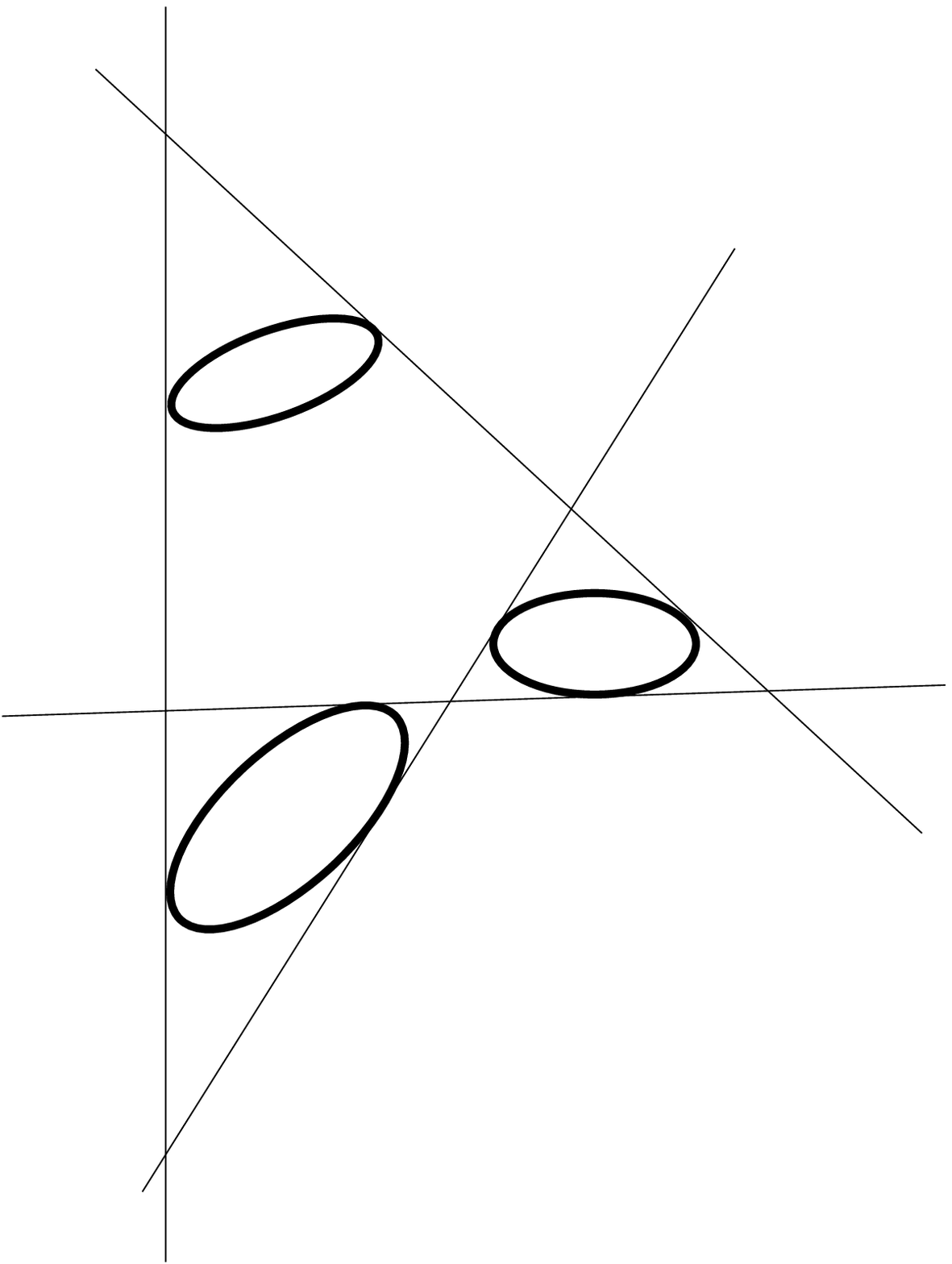}
  \caption{Hessian and flecnodal cur\-ves for a factorisable polynomial.}
  \label{factorisableafin} 
  \end{minipage}
  \hskip 0.7cm
  \begin{minipage}{0.47\textwidth}
    \includegraphics[width=2.0in]{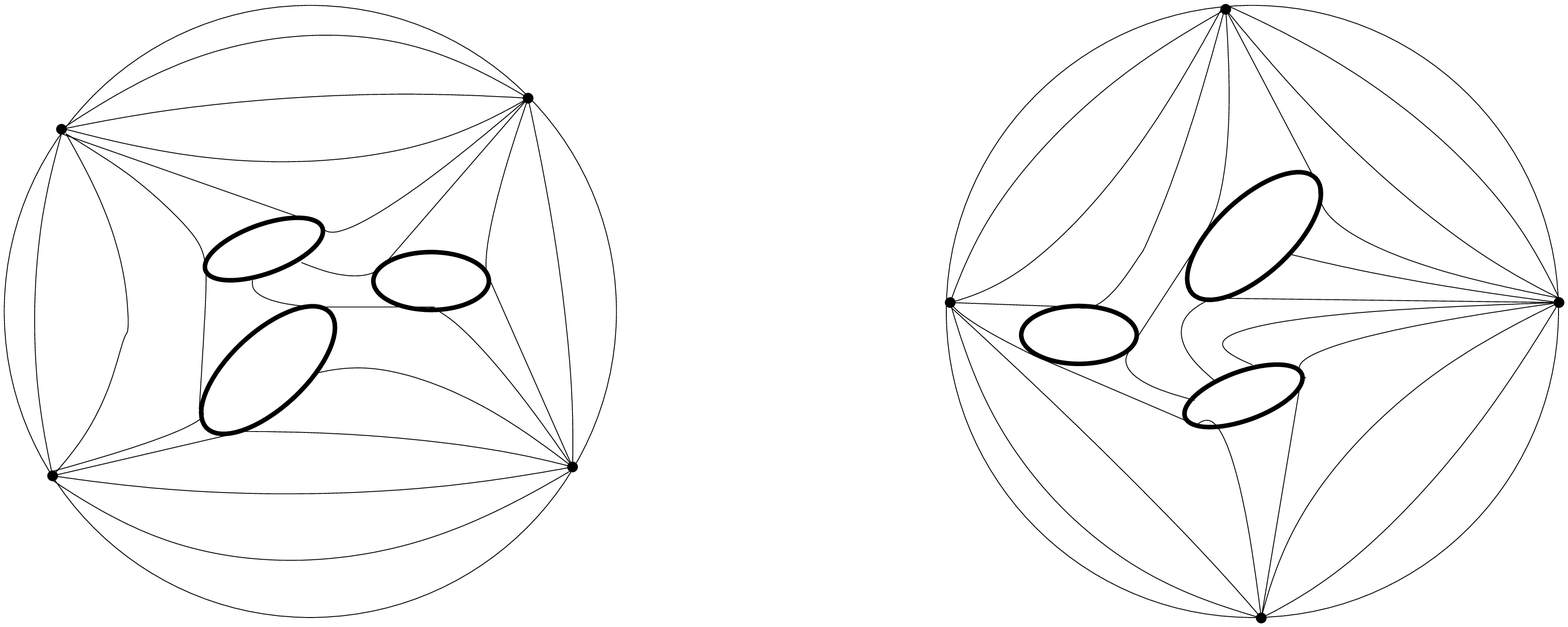}
    \caption{Foliation of $\mathbb{Y}_{k}$ for a quartic factorisable polynomial.}
    \label{ejemplog4}
  \end{minipage}
\end{figure}

\noindent In Fig. \ref{factorisableafin}, we show the affine geometrical 
structure of the quartic factorisable polynomial 
$g(x,y) = y (x+3) (x-y) (y+x-3)$. 
By Remark \ref{indexatprojplane}, each field $\widetilde{\mathbb{X}_k}$ 
associated to this example has two singular points at infinity while each 
$\mathbb{Y}_k$ has eight singular points at infinity.
We conclude this example by offering in Fig. \ref{ejemplog4} a picture of the 
qualitative behaviour of $\mathbb{Y}_{k}$ in the two hemispheres.}
\end{examples}


\section{Upper bounds for the number of godrons}\label{sec:6}
 
In the first part of this chapter, we prove a formula that relates the Euler 
characteristic of the surface $B^{\pm}$ with the Poincar\'e indices at 
singularities of $\widetilde{\mathbb{X}}$ when such field is defined on 
$B^{\pm}$, respectively. In the second part, as an application, we give an upper 
bound for the number of godrons when the projective Hessian curve of $f$
is only constituted by exterior ovals. Before stating our results, we introduce 
some definitions.

The tangency of the asymptotic line with the Hessian curve of $f$ at a 
godron is either, {\it exterior} or {\it interior} 
\cite{Bancthom}. In the first case, we say that such godron has an 
\emph{interior tangency} and in the second case, an \emph{exterior tangency.}

\begin{theorem}\label{thglobalbminus}
Let $f$ $\in
\mathbb{R} \left[ x,y\right] $ be a polynomial of degree $n\,$ whose graph $S_f$ 
is generic and its 
projective Hessian curve is not tangent to the line at infinity. 
Assume that $f_{n}$ has no repeated factors and that $\widetilde{\mathbb{X}}$ 
is a projective extension of a field of asymptotic directions. Then
\begin{eqnarray*}
\sum_{\xi \in \Sing (\widetilde{\mathbb{X}})} \Ind_\xi\left( 
\widetilde{\mathbb{X}} \right) =
\chi \left( B^{\epsilon}\right) +\frac{%
	P_{i}-P_{e}}{2},
\end{eqnarray*} 
where $\epsilon$ is either $+$ or $-$ and $\widetilde{\mathbb{X}}$ is defined on 
$B^{\epsilon}$. 
In both cases, $P_{i}$ denotes the number of godrons with an interior
tangency and $P_{e}$, with an exterior tangency.
\end{theorem}

\noindent{\bf Proof.} Since the projective Hessian curve is transversal to the 
line at infinity $z=0$, all tangencies of $\widetilde{\mathbb{X}}$ with the 
projective Hessian curve occur in the Hessian curve of $f$, Hess $f(x,y)=0$.

Suppose that $\widetilde{\mathbb{X}}$ is defined on the smooth surface $B^{-}$. 
This surface 
is composed by a finite number of orientable connected components denoted by 
$\,D_{1},\ldots , D_{s}$ and a connected component $\mathbb{M}$  homeomorphic to a 
closed M\"{o}bius strip with a finite number of open discs removed. 
For $l=1, \ldots ,s$, we denote by $P_{i}^{D_{l}}$ and $P_{e}^{D_{l}}$
the number of godrons having an interior and exterior tangency on the boundary of
$\,D_{l}$, respectively. Poincar\'e-Hopf's Theorem for surfaces with boundary implies
\begin{eqnarray}\label{PHorientable}
\sum_{l=1}^s \,\,\left(\sum_{\substack{\xi\in \Sing (\widetilde{\mathbb{X}})  \\ 
\xi\in D_l\mbox{\qquad}}} \Ind_\xi (\widetilde{\mathbb{X}})\right) =  \sum_{l=1}^s 
\chi \left(D_l\right) + \sum_{l=1}^s \frac{P_{i}^{D_{l}}-P_{e}^{D_{l}}}{2}.
\end{eqnarray}

Now, we shall prove a version of Poincar\'e-Hopf's 
Theorem for the nonorientable surface $\mathbb{M}$. 
The projective extension $\widetilde{\mathbb{X}}$ is the restriction of a field 
$\mathbb{Y}$ defined by (\ref{EDLA}) to a hemisphere. Such field $\mathbb{Y}$ is 
defined on
an orientable surface $\mathbb{DM} \subset \mathbb{S}^2$ which is a double covering 
of $\mathbb{M}$. So, $\, \chi \left( \mathbb{DM}\right) = 2 \chi 
(\mathbb{M})$. By considering the Poincar\'e-Hopf Theorem for the field $\mathbb{Y}$,
\begin{equation}\label{PHdcovering}
\sum_{\substack{\xi \in \Sing (\mathbb{Y})\\ \xi \in \mathbb{DM}
\mbox{\qquad}}} \Ind_\xi \left( \mathbb{Y} \right) =\chi \left( \mathbb{DM}\right) 
+\frac{P_{i}^{\mathbb{DM}}- P_e^{\mathbb{DM}}}{2}.
\end{equation}
Since the number of tangencies that the field $\mathbb{Y}$ has with the boundary 
of $\mathbb{DM}$ is twice the number of tangencies of $\widetilde{\mathbb{X}}$ 
with the boundary of $\mathbb{M}$, we obtain the relation
 $\,P_{i}^{\mathbb{DM}}- P_e^{\mathbb{DM}} = 2 (P_{i}^{\mathbb{M}}- 
P_e^{\mathbb{M}}).$ Moreover, by Remark \ref{foliacionesenesfera}, 
$\,\sum_{\xi \in \Sing (\mathbb{Y})} \Ind_\xi \left( \mathbb{Y} \right) = 2\sum_{\xi \in
 \Sing (\widetilde{\mathbb{X}})} \Ind_\xi ( \widetilde{\mathbb{X}}).$
By replacing these expressions in (\ref{PHdcovering})
we obtain the desired equality 
\begin{equation}\label{poincare-moebius-band}
\sum_{\substack{\xi \in \Sing (\widetilde{\mathbb{X}})\\ \xi 
		\in \mathbb{M}\mbox{\qquad}}} \Ind_\xi\left(\widetilde{\mathbb{X}}
\right) =\chi \left( \mathbb{M}\right)
+\frac{P_{i}^{\mathbb{M}}-P_{e}^{\mathbb{M}}}{2}.
\end{equation}

\noindent The proof concludes from (\ref{PHorientable}) and 
(\ref{poincare-moebius-band}) since $\,B^{-}=\mathbb{M}\amalg D_{1}\amalg 
\cdots \amalg D_{s}$.
\hfill $\Box$
\medskip

Points on a generic algebraic surface in $\mathbb{CP}^{3}$ are also classified
in terms of the maximum order of contact of the tangent lines at them with the 
surface. George Salmon proves in \cite{salmon} that such a surface of 
degree $n$ has $\, 2n(n-2)(11n-24)\,$ godrons (points at which the asymptotic line is
tangent to the parabolic curve \cite{McCry&SHrfn}). This number is an upper 
bound for the number of godrons on a generic algebraic surface in $
\mathbb{RP}^{3}$. When the graph of a polynomial $f\in\mathbb{R}
\left[ x,y\right]$ of degree $n$ is generic, an upper bound for the number of godrons
is given (\cite{H&O&S1}, Theorem 5), namely,
\begin{equation}\label{csnppe}
\#\, \{ \mbox{Godrons in } S_f \} \leq \left( n-2\right) \left( 5n-12\right).
\end{equation}

\begin{corollary}\label{pssuperiorbounds}
Let $f$ $\in\mathbb{R} \left[ x,y\right] $ be a polynomial of degree $n$ such 
that $S_f$ is generic. Suppose that the projective
Hessian curve is not tangent to the line at infinity. If the polynomial $f_n$ 
has $k$ distinct real linear factors, then
$$P_{i}\leq \frac{( n-2)(8n-21) +k}{2} \qquad \mbox{and} \qquad	
P_{e}\leq  1 + \frac{(n-2)(8n-21)- k}{2}.$$
\end{corollary}

\noindent {\bf Proof.} Since $S_f$ is generic, the projective Hessian curve of $f$
is an algebraic curve of degree $2n-4$. By Theorem \ref{petrows}, the value 
$\chi (B^{+})$ satisfies: 
$$-\frac{3(n-2)(n-3)}{2} \leq \chi \left( B^+\right) \leq 1+ \frac{3(n-2)(n-3)}{2}.$$

\noindent Because 
$\chi ( B^{+}) =1-\chi ( B^{-})$, we have that $\chi (B^{-})$ satisfies the 
inequalities $\, -\frac{3}{2}(n-2)(n-3) \leq \chi ( B^-) \leq 1 +
\frac{3}{2}(n-2)(n-3).$ In conclusion, we obtain
\begin{equation}\label{inequalityderpetrovsjy} 
- \frac{3(n-2)(n-3)}{2} -1 \leq -\chi \left( B^{\pm}\right) \leq \frac{3(n-2)(n-3)}
{2}.
\end{equation}

\noindent If the set $b^{\pm}\cap\mathbb{R}^{2}$ 
is hyperbolic, then, by Theorem \ref{thglobalbminus} 
\begin{equation}\label{pipet4}
\frac{P_{i}-P_{e}}{2} = \sum_{\xi \in \Sing
	(\widetilde{\mathbb{X}})} \Ind_\xi\left(\widetilde{\mathbb{X}} \right) -
\chi \left( B^{\pm}\right).
\end{equation}

\noindent According to Corollary \ref{lemma01} and Theorem \ref{thindexonehalf}, 
$\, \sum_{\xi \in \Sing (\widetilde{\mathbb{X}})} 
\Ind_\xi (\widetilde{\mathbb{X}}) = \frac{k}{2}$.
Therefore, by adding $\frac{k}{2}$ to the inequalities  
(\ref{inequalityderpetrovsjy}) and using (\ref{pipet4}), we get
\begin{equation}\label{pipe}
- 3 \left( n-2\right) \left( n-3\right) -2 +k \leq P_{i}-P_{e}\leq 3 
\left( n-2\right) \left( n-3\right) +k.
\end{equation} 
The proof follows from inequalities (\ref{csnppe}) and (\ref{pipe}).
\hfill $\Box$
\medskip

When the Hessian curve of $f$ is a convex compact curve and the set
 $b^-\cap\mathbb{R}^{2}$
is hyperbolic, the second author of this paper joined to L.I. Hern\'andez-Mart\'\i 
nez and F. S\'anchez-Bringas to prove that $\,n(3n-14)+18\,$ is an upper bound 
for the number of godrons lying on the boundary of the unbounded connected 
component $C_u$ (\cite{H&O&S2}, Theorem 10). In the 
following result, we improve such bound under different assumptions
and we analyse the unbounded case: we give an upper bound for the 
number of godrons that are on the boundary of $\mathbb{M}$.

\begin{theorem}\label{PSP2014} 
Let $f\in \mathbb{R}
\left[ x,y\right] $ be a polynomial of degree $n$ whose graph $S_f$ is generic. 
Suppose that the projective Hessian curve of $f$, constituted only by exterior 
ovals, is convex and it is not tangent to the line at infinity. If $H$ is 
contained in $B^-$ and the polynomial $f_n$ has $k$ distinct real linear 
factors, then the maximal number of godrons is $\, 3 (n-2)(n-3) + k.$
\end{theorem}

\noindent {\bf Proof.}
On the one hand, since any projective extension $\widetilde{\mathbb{X}}$ of a 
field of asymptotic lines is defined on $B^-$ the
expression $P_i-P_e$ satisfies the second inequality of (\ref{pipe}), that is,
$\, P_{i}-P_{e}\leq 3 \left( n-2\right) \left( n-3\right) +k.$
On the other hand, all godrons have an interior tangency because the projective 
Hessian curve of $f$ is convex and the set $b^{-}\cap \mathbb{R}^{2}$ is hyperbolic.
Therefore, $P_{e} = 0$ and  $\,P_i $ equals the total number of godrons. 
\hfill$\Box$
\medskip
\section{Appendix}\label{apend}
\noindent {\bf Proof of Theorem \protect\ref{thindexonehalf}}.
Let $p \in \mathbb{S}^2 \subset \mathbb{R}^3=\{(u,v,\omega)\}$ be a
singular point at infinity of the field $\, \mathbb{Y}_k, k=1,2$. According to 
Corollary \ref{lemma01},  a real linear factor $l$ of $f_n$ defines the point $p$
and by hypothesis, the multiplicity of $l$ is one.
After a suitable linear change of coordinates on
the $xy$-plane, we have that $l(x,y) = y$, $\, p=(1,0,0)$ and 
\begin{eqnarray}\label{desarrollofn}
f_{n}\left( x,y\right) &=&y\left( \sum_{\substack{ i=0 }}^{n-1}a_{i,n-i} \,
x^{i}y^{n-1-i}\right),\text{ with }a_{n-1,1}\neq 0.
	\label{polhomanminus1} 
\end{eqnarray}
In the chart $u=1$ the fields $\mathbb{Y}_{k}, k=1,2,$ restricted to the set 
$\{(u,v,\omega)\in \mathbb{S}^2 | u> 0\}$ are described by the quadratic equation
\begin{equation}\label{formaencartaafin}
\omega ^{2}F_{vv}( 1,v,\omega)(dv)^2 + 2 \omega B\left(1,v,\omega \right) dv d\omega
+ S\left( 1,v,\omega \right)(d\omega)^2 = 0.	
\end{equation}%

\noindent By Theorem \ref{puntos-sing-infinito}, $S(1,0,0) = 0$ and $H_f(p) <0$.
Moreover, in a neighbourhood of the origin on the $v\omega$-plane the two fields of directions defined by (\ref{formaencartaafin}) are described by (see equation (\ref{campos-separados}))
\begin{eqnarray}\label{dfields2}
	R_{k}\left( v,\omega \right) dv + 2 S\left(1, v,\omega \right) d\omega =
	0,\,\,\, k=1,2.
\end{eqnarray}

\noindent We denote by $\mathcal{G}_{1}$ and $\mathcal{G}_{2}$ the foliations of 
these fields. The proof is based on the following geometric idea. 
Consider the sets
$W_{U}=\left\{
\left. \left( v,\omega \right) \in \mathbb{R}^{2}
\right\vert \omega > 0\right\} $ and 
$W_{L}=\left\{ \left. \left( v,\omega \right) \in 
\mathbb{R}^{2}\right\vert \omega <
0\right\} $. The key point is to prove that there exists a
neighbourhood of the origin, denoted by $W \subset\mathbb{R}^{2}$,  
at which one of the two foliations, $\mathcal{G}_{1}$ for example, is tangent, in 
$W\cap W_{U}$, to a vector 
field having a node at the the origin and, in $W\cap W_{L}$, is tangent to a 
nonsingular vector field. Simultaneously, we will have that the foliation
$\mathcal{G}_{2}$ is tangent to the same vector fields, but in this case in 
the sets $W\cap W_{L}$ and $W\cap W_{U}$, respectively.

From the expression of the fields described in (\ref{dfields2}) and setting
$\tilde{S}(v,\omega) = S\left( 1,v,\omega \right )$,
we define the following vector fields on the $v\omega$-plane which have similar 
qualitative behaviours.
\[
Y_{k}\left( v,\omega \right) =\left(\tilde{S}\left( v,\omega \right) ,\,\omega\,
T_{k}\left( v,\omega \right) \right), \,\,\,\, k=1,2,
\]%
where $\,
T_{k}\left( v,\omega \right) =-2B\left( 1,v,\omega \right) + 2
\left( -1\right) ^{k}\sqrt{- H_f(1,v,\omega)}.$

It is clear that in a punctured neighbourhood of the origin the foliation 
$\mathcal{G}_{1}$ is tangent to the vector field $Y_{1}$ if $\omega > 0,$ and 
tangent to the vector field $
Y_{2}$ if $\omega < 0.$ Respectively, the foliation $\mathcal{G}_{2}$ is
tangent to the vector field $Y_{2}$ when $\omega > 0,$ and it is tangent to
the vector field $Y_{1}$ for $\omega < 0.$

Now, we shall describe the topological type of the origin. A straightforward 
calculation shows the equality

\begin{equation}\label{relacionT1T2}
T_{1}\left( v,\omega \right) T_{2}\left( v,\omega \right) =4\,
\tilde{S}\left( v,\omega \right)  F_{vv}\left(1,v, \omega\right).
\end{equation}

\noindent Because
\begin{eqnarray*}
&&B(1,0,0) = -\frac{\partial^2 f_n}{\partial u \partial v}|_{(1,0)} = - 
(n-1)a_{n-1,1}\quad \mbox{ and } \\ 
&&H_f(1,0,0) =  Hessf_{n}(1,0) = -\left( n-1\right)^2 a_{n-1,1}^{2},
\end{eqnarray*}
we assert that 
$$ T_{k}\left( 0,0\right) = 2\left( n-1\right)a_{n-1,1} + 2\left( -1\right)^{k} 
\left(n-1\right)\left\vert a_{n-1,1}\right\vert .$$
Therefore, if $a_{n-1,1}>0,$ then $T_{1}\left( 0,0\right) =0$ and $
T_{2}\left( 0,0\right) =4\left( n-1\right) a_{n-1,1}$. In case $a_{n-1,1}<0,$
$T_{2}\left( 0,0\right) =0$ and $T_{1}\left( 0,0\right) =4\left(n-1\right) a_{n-1,1}.$
It allows us to analyse the linear part of the vector field $Y_{k} $ 
at the origin. 

\noindent $DY_{k}\left( 0,0\right) =$
\begin{align*}
&=\left. \left( 
\begin{array}{cc}
2\frac{\partial }{\partial v}\tilde{S}\left( v,\omega \right)  
&2\frac{\partial }{%
	\partial \omega }\tilde{S}\left( v,\omega \right)  \\ 
\omega \frac{\partial }{\partial v} T_{k}\left( v,\omega \right)& 
\omega \frac{\partial }{\partial \omega} T_k(v,\omega) + T_{k}\left( v,
\omega \right) 
\end{array}%
\right) \right\vert _{\left( 0,0\right) } \\
&=\left. \left( 
\begin{array}{cc}
2 n(n-1) \frac{\partial}{\partial v} \sum_{i=0}^n \omega^{n-i} f_{i}\left(
1,v\right)  & 2 n(n-1) \frac{\partial}{\partial \omega} \sum_{i=0}^n 
\omega^{n-i} f_{i}(1,v)\\ 
\omega \frac{\partial}{\partial v}T_{k}(v,\omega)
& \omega \frac{\partial}{\partial \omega} T_k(v,\omega) + T_{k}( 
v,\omega) 
\end{array}%
\right) \right\vert _{\left( 0,0\right) } \\
&=\left( 
\begin{array}{cc}
2n\left( n-1\right) a_{n-1,1} & 2 n (n-1)\frac{\partial }{\partial 
	\omega } f_{n-1}\left( 1,v\right)\vert_{v=0} \\ 
0 & 2\left( n-1\right)a_{n-1,1} +  2\left( n-1\right) (-1)^k \left\vert 
a_{n-1,1}\right\vert 
\end{array}%
\right).
\end{align*}

\noindent The matrix $DY_{1}\left( 0,0\right) $ ($DY_{2}\left( 0,0\right)$)
has two nonzero real eigenvalues with the same sign if $a_{n-1,1} <0$ (if 
$a_{n-1,1} >0$ respectively). In conclusion, if $\,a_{n-1,1}\,$ is positive 
(negative) then the origin is a singular point of type node of the vector field 
$Y_{2}\, ($respectively $Y_{1})$.

Now, suppose $a_{n-1,1}>0$ (the negative case is analogous) and consider the vector 
field
\[
Z_{1}\left( v,\omega \right) =\left( T_{2}\left( v,\omega \right) , 4\, \omega
F_{vv}\left( 1,v, \omega\right) \right). 
\]%

\noindent Since $T_{2}\left( 0,0\right) \neq 0$, the origin is a nonsingular point of 
this vector field. Moreover, by (\ref{relacionT1T2}) this field satisfies the equality
\[
T_{2}\left( v,\omega \right) \, Y_{1}\left( v,\omega \right) = \tilde{S} \left(
v,\omega \right) \, Z_{1}\left( v,\omega \right).
\]
This relation implies that $Z_{1}$ is tangent to the
foliation of $Y_{1}$.
\hfill$\Box$ 
\medskip

\noindent {\bf Acknowledgments}
This research was partially supported by DGAPA-UNAM grant PAPIIT-IN111415 and by 
CONACyT grant CB$\_$219722. We are truly greatful with the referee for the 
suggested changes to this work.

\end{document}